\documentclass[11pt]{article}
\usepackage{graphicx,graphics}
\usepackage{xcolor}
\usepackage{lineno}
\usepackage[unicode=true,pdfusetitle,
bookmarks=true,bookmarksnumbered=false,bookmarksopen=false,
breaklinks=false,pdfborder={0 0 0},backref=false,colorlinks=false]
{hyperref}

\usepackage{geometry}
\usepackage[english]{babel}
\usepackage{amsmath}
\usepackage{amsthm}
\usepackage{amssymb}
\usepackage{amsfonts}
\usepackage{graphicx}
\usepackage{hyperref}
\usepackage{tikz}
\usepackage{tikz-3dplot}
\usepackage{pictexwd, dcpic}
\usepackage{pgfplots}
\pgfplotsset{compat=newest}
\usepackage{caption}
\usepackage{subcaption}
\usepackage{bigints}
\usepackage{tipa}
\usepackage{lipsum}

\newcommand\blfootnote[1]{%
  \begingroup
  \renewcommand\thefootnote{}\footnote{#1}%
  \addtocounter{footnote}{-1}%
  \endgroup
}
\newtheorem{theorem}{Theorem}[section]

\newtheorem{prop}[theorem]{Proposition}

\newtheorem{remark}{Remark}[section]

\begin{document}
%\begin{linenumbers}
\begin{center}  {\Large \bf An integral renewal equation approach to behavioural epidemic models with information index}
\end{center}
\smallskip
\begin{center}
{\small \textsc{Bruno Buonomo$^{1,*}$}, 
\textsc{Eleonora Messina$^{1,*}$},
\textsc{Claudia Panico$^{1,*}$},
\textsc{Antonia Vecchio$^{2,*}$}}
\end{center}

\blfootnote{$^* $ \textit{ E-mail addresses:} buonomo@unina.it (B.Buonomo), eleonora.messina@unina.it (E. Messina), claudia.panico2@unina.it (C. Panico), antonia.vecchio@cnr.it (A. Vecchio).}

\begin{center} {\small \sl $^{1}$ Department of Mathematics and Applications,
University of Naples Federico II,\\ via Cintia, I-80126
Naples Italy.}\\
{\small \sl $^2$National Research Council of Italy, Institute for Computational Application ``M. Picone'',\\ Via P. Castellino, 111 - I-80131 Naples, Italy.}\\
%{\small \sl $*$ Member of the INdAM Research Group GNCS, Italy}
\end{center}

\medskip

{\centerline{\bf Abstract}
	\begin{quote}
		
		{\small We propose an integral model describing an epidemic of an infectious disease. The model is behavioural in the sense that the force of infection includes a distributed delay, called {\it information index}, that describes the opinion-driven human behavioural changes. The information index, in turn, contains a memory kernel to mimic how the individuals maintain memory of the past values of the infection. We obtain sufficient conditions for the endemic equilibrium to be locally stable. In particular, we show that when the infectivity function is represented by an exponential distribution, stability is guaranteed  by the weak Erlang memory kernel. However, through numerical simulations, we show that oscillations, possibly self-sustained,  may arise when the memory is more focused in the disease's past history, as exemplified by the strong Erlang kernel. We also show the model solutions in cases of different infectivity functions taken from studies where specific diseases like Influenza and SARS are considered.
		}
\end{quote}}

%\medskip
%
%\noindent { {\bf Subject class}: }
%
\medskip


\normalsize

%\maketitle

\section{Introduction}

Models of Mathematical Epidemiology (ME) have proved to be an effective tool for understanding epidemic dynamics and for identifying the best intervention and mitigation strategies. Many noteworthy cases can be cited where models of ME have been applied to recent diseases and parametrized with official data. Clear examples are the recent COVID-19 pandemic \cite{dellaRossa2020,Dorigatti2020,Gatto2020,Pellis_2022}, the transmission of Middle East Respiratory Syndrome (MERS) in South Korea and Saudi Arabia \cite{RUI2022102243}, and the spread of the West Nile virus in Italy \cite{fesce2023}.

In all of the above-mentioned cases,  \textit{compartmental} models were used in which the individuals of a given population were divided in mutually exclusive groups (the \textit{compartments}) according to their respective disease status. Although compartmental models applied to real epidemics often have a complex structure, they are in fact extensions of {\it minimal} epidemic models. By `minimal', we mean models with a reduced mathematical complexity where the epidemic is described through slight variants of the basic SIR (Susceptible–Infectious–Removed) or SEIR (Susceptible–Exposed–Infectious–Removed) models. 
 Study of these minimal models have attracted the attention of researchers for many years and may prove critical as paradigmatic models of specific dynamical outcomes such as bifurcations, hysteresis and oscillations \cite{CAPASSO1993,CAPASSO197843,Martcheva_2015}. It is well known that the minimal compartmental models refer back to the work of Kermack and McKendrick. In 1927, they published a celebrated paper in which an integral \textit{renewal equation} (i.e. a model for epidemic propagation based on the individual infectiousness) was introduced \cite{KM1927} (see also the reprint \cite{KERMACK199133}). %Their model has been reformulated later, in a modern fashion, in terms of the \textit{constitutive equation} for the Force of Infection (FoI), i.e. the rate at which a susceptible becomes infected \cite{JBD2012}. 
 Many researchers have made substantial contributions to the development of integral epidemic models in the field of epidemiology. An in-depth coverage of integral models in the mathematical modeling of infectious diseases can be found in current classical books on ME \cite{Brauer_Castillo-Chavez_Feng_2019, DIEKHEESTE}.
 Applications of integral epidemic models to actual diseases, such as influenza and Severe Acute Respiratory Syndrome (SARS), were provided by M. G. Roberts and coworkers \cite{Roberts2006,Roberts2007}. 
 O. Diekmann and coworkers proposed mathematical formulations for a wide range of epidemic models that incorporate demography and non-permanent immunity as a scalar renewal equation for the Force of Infection (FoI), i.e. the rate at which a susceptible becomes infected \cite{JBD2012}. In such a paper, the concept is introduced  of \textit{infectivity function}, which represents the expected contribution to the FoI by an individual infected since $\tau$ units of time. Kermack and McKendrick’s integral renewal equation reduces to classical models, such as SIR and SEIR, when the infectivity function assumes specific forms \cite{Finite2020,DIEKIN}.
 
Both renewal equations and  compartmental models of ME do exhibit a notable limitation, however, that in the context of  contemporary societies cannot be disregarded: they do not take into account the feedback effects produced by the behavioural change of individuals during an epidemic. Classical compartmental models have been revisited in recent  years within the framework of a new field called \textit{Behavioural Epidemiology of Infectious Disease} (BEID) \cite{mado2013,WANG20161}. Among the possible approaches within the BEID, the one based on the so-called \textit{information index} (i.e.  a distributed delay which describes the opinion-driven changes of behaviour of individuals) is particularly promising and increasingly used \cite{Bai2014, KUMAR201977, LopezCruz2022}.
The first studies on minimal epidemic models embedding the information index date back to 2007 with A. d'Onofrio and coworkers \cite{donofrio2007}. Generally speaking, scholars working on models of BEID focus on potential complex dynamics caused by the presence of behavioural effects. For example, it has been proven that the opinion-driven change of behaviour, modeled through the information index, can induce oscillatory dynamics that go uncaptured by analogous non-behavioural models \cite{donofrio2009,donofrio2007,mado2013,WANG20161}. A key factor in the occurrence of oscillations is the type of \textit{memory kernel} (i.e. the function used to model how a population keeps memory of information and rumours regarding a disease \cite{donofrio2009,DONOFRIO2022,donofrio2007}). Recently, minimal epidemic models, including the information index, have been extended to incorporate more complex models for the study of vaccine hesitancy \cite{buonomo2022} and to obtain the value of the information-related parameters by using official data released by public health authorities  \cite{BBRDMCovid,BBRDMeningitis}. It is worth mentioning that in several papers, an \textit{awareness variable} is employed to represent behavioural responses due to information. Such a variable can be considered as a special case of the information index \cite{LACITIGNOLA2021110739, Saha2020, ZUO2023127905}. 

 In this paper, motivated by the recent promising advances of BEID, we offer a new contribution to the theory of integral renewal equations in ME. We propose a \textit{behavioural} renewal equation where the FoI incorporates the information index and the model therefore takes into account human behavioural feedback. As far as we know, this is the first investigation on integral renewal equations that embed the information index.
 
 Our approach provides a qualitative analysis of the model based on stability theory. That is, we obtain sufficient conditions for the endemic equilibrium to be locally stable. Such conditions are written in terms of the memory kernel, the infectivity function and the so-called \textit{message function}, which describes the information that forms the opinion of individuals. Moreover,  through numerical simulations based on a non-standard method, we show how the memory kernel affects the dynamical outcomes of the model. We also show the model solutions when the infectivity function has the form used to model the spread of infectious diseases like 
 %in case of some different infectivity functions describing infectious diseases, like
 Influenza and SARS. Since the renewal equation is a generalization of epidemic models ruled by Ordinary Differential Equations (ODEs), this study also extends the results concerning the minimal behavioural models to more general (and potentially more realistic) cases.

The paper is organized as follows: in Section \ref{sec2}, we build the renewal equation model in detail. Section \ref{sec3} focuses on the basic properties of the solution, including positivity and boundedness. We prove that the system admits two equilibria under specific assumptions.  In Section \ref{sec4}, we derive the characteristic equation as a preliminary for the stability analysis. Section \ref{sec5} is devoted to the local stability analysis of the equilibria. In Section \ref{sec6}, we shortly discuss the non-standard method built for the problem at hand and show some applications using infectivity functions proposed for specific infectious diseases. Finally, conclusions and future prospectives are provided in Section \ref{sec7}.

\section{The behavioural renewal equation}\label{sec2}
In the formulation of the Kermack and McKendrick's integral renewal equation \cite{JBD2012}, the FoI at time $t$, denoted by $\hat F(t)$, is introduced through the following constitutive equation:

\begin{equation}\label{ConstEqHat}
 \hat F(t) = \int_0^{\infty}  S(t-\tau) \hat F(t-\tau) A_{\mu}(\tau)  d \tau.
 \end{equation}

\noindent Here, $S$ denotes the size of the susceptible compartment;  $A_{\mu}(\tau)= e^{-\mu \tau} A(\tau)$, where $A(\tau)$ is the \textit{Infectivity Function} (IF), that is the expected contribution to the FoI by an individual that was itself infected $\tau$ units of time ago; $\mu$ is the per capita natural death rate, and $e^{-\mu \tau}$ represents the probability that an individual stays alive for at least $\tau$ units of time \cite{BRAUER2016}. We call $A_{\mu}(\tau)$ the \textit{Infectivity Function with Demography} (IFD).  The IF is assumed to be an $L^1$ function from $[0, \infty)$ to $[0, \infty)$.

The balance equation describing the changes of the variable $S$, when demographic turnover is taking into account, is given by:
\begin{equation}
\label{EqSnonb}
\dot S (t) = \lambda - \mu S (t) -  \hat F(t) S(t),
\end{equation}
where the upper dot denotes the time derivative, $\lambda$ represents the inflow of susceptibles, $\mu$ is the natural death rate and the quantity $\hat F(t) S$ is the \textit{incidence}, that is the number of new cases per unit time and area \cite{JBD2012}. Therefore, the constitutive equation \eqref{ConstEqHat}  describes how the current FoI depends on past incidence.

It is well known that compartmental models can be derived from the renewal equation \eqref{ConstEqHat} (see \cite{Finite2020,DIEKIN}). In this regard, a key role is played by the IF. For example, if $A(\tau)$ is an exponential function or a linear combination of exponential functions, from \eqref{ConstEqHat} one can derive the classical SIR and SEIR compartmental models ruled by three and four ODEs, respectively. 

In the behavioural formulation of epidemic models based on the information index it is explicitly taken into account that the transmission of information and rumours about the disease are often preceded by articulated routine procedures so that public awareness regarding the disease takes time to build \cite{mado2013, WANG20161}. Moreover, people keep memory of the past values of the infection. Therefore, delay and memory are two major aspects of the opinion-driven change of behaviour. The information index, $M(t)$, summarizes these aspects and is defined as the following distributed delay:

	\begin{equation}\label{emme}
    M(t)= \int_{0}^{\infty}  H(t-\tau) \, g \big( x_1(t-\tau), x_2(t-\tau), \dots, x_n(t-\tau) \big)  \, K(\tau)\, d \tau,
\end{equation}
where $H(\cdot)$ is the Heaviside step function and 
\begin{itemize}
    \item[\textit{(i)}] the quantities $x_1, x_2, \dots, x_n$ are a subset of the state variables of the epidemic model under consideration;
    
     \item[\textit{(ii)}] the message function $g$ describes the individuals' perceived risks associated to the disease. In principle, it may depend on all the state variables, but, in practice, it is assumed to depend only on some specific quantities as, for example, the incidence and the prevalence of the disease; 
 
     \item[\textit{(iii)}]  $K(\cdot)$ is the memory kernel, such that $K(t) >0$ and $\int_0^{\infty} K(t) dt=1$. For example, the memory kernel can be chosen to be a Dirac delta $K(t)=\delta(t)$, which represents the instantaneous propagation of the information, or an Erlang distribution with rate parameter $a$ and shape parameter $n$,
     
     $$
     K(t)=Erl_{n,a}(t)= \frac{a^n}{(n-1)!} t^{n-1} e^{-n t}, \quad a \in \mathbb{R^+}, \quad n \in \mathbb{N^+},
     $$
     which includes both the cases of unimodal kernels (when $n>1$) and of exponentially fading memory (when $n=1$) in which the memory of the mean individual decays with characteristic length $T = 1/a$, see \cite{mado2013,WANG20161}.
 \end{itemize}

\noindent Now, going back to the balance equation (\ref{EqSnonb}), we assume that the FoI, $\hat F$, is made up of two terms:
$$
	\hat F (t) = \beta (M(t))\, F(t).
$$
The first term, $\beta(M(t))$, is a function that describes the inhibitory contribution to the FoI due to the information and rumours  regarding the disease status. The function $\beta$ is assumed to verify the following conditions:  $\beta(M)>0$, $\beta(0)=1$ and $\beta^{\prime}(M) <0$. The second term, $F(t)$, is the \textit{driving factor of the FoI} which we introduce through the following constitutive law:
\begin{equation}\label{BeConstEq}
	F(t) = \int_0^{\infty}  \beta \big(M(t-\tau)\big) S(t-\tau) F(t-\tau) A_{\mu}(\tau) d \tau.
\end{equation}
The balance equation for the susceptibles is still equation (\ref{EqSnonb}), which can now be written:
\begin{equation}\label{EqS}
	\dot S(t) =\lambda - \mu S(t) -  \beta(M(t)) F(t) S(t).
\end{equation}
Note that the formulation (\ref{BeConstEq}) is a ``behavioural'' variant of (\ref{ConstEqHat}) and tells us how the driving factor of the FoI depends on the past incidence. We underline again that $F$ is not a FoI. However, it becomes  a FoI in the ``non-- behavioural'' case $M=0$.

It remains now to discuss the  message function $g$ in (\ref{emme}).  We propose that it depends on the quantity $F$. In order words, we assume that the perception of individuals regarding the disease is based on the awareness regarding the driving factor of the FoI. This assumption is based on the idea that the quantity $F$ generalizes in some extent the prevalence of the disease (i.e. the size of infectious population, say $I$) as shown for a specific example in Proposition \ref{prop:1} below. Therefore,  the assumption that $g$ depends on $F$ generalizes the idea of prevalence-dependent behavioural epidemic models studied through ODEs in \cite{donofrio2009,donofrio2020}.

The discussion above leads to the following form of the information index:
 \begin{equation}
 	\label{BeM}
      M(t)= \int_{0}^{\infty}  H(t-\tau)  g \big(F(t-\tau)\big)\,K(\tau)\,d \tau,
 \end{equation}

\noindent where  $g(0)=0$, $g(F)\geq 0$ and $g^{\prime}(F) >0$.  Collecting together  \eqref{BeConstEq}, \eqref{BeM} and \eqref{EqS}, the complete model reads

 \begin{equation}\label{mod continuo22}
     \begin{aligned}
         &\dot S(t)=\lambda - \mu S(t) - \beta(M(t)) S(t) F(t), \\
         &F(t) = \int_{0}^{\infty} \beta \big(M(t-\tau)\big) S(t-\tau) F(t-\tau)  A_{\mu}(\tau)   d \tau, \\
          &M(t)= \int_{0}^{\infty}   H(t-\tau) g \big(F(t-\tau)\big)K(\tau)  d \tau.
     \end{aligned}
 \end{equation}

\noindent For the sake of clarity, we report in Table \ref{tab1} the definition of the parameters and quantities related to model \eqref{mod continuo22}  and in Table \ref{tab:assump} the assumptions on the functions involved in model \eqref{mod continuo22}, which we suppose to hold from now on.

\begin{table}[t]
\centering
\small
\begin{tabular}{|c |l |}\hline
 Symbol	 &  Meaning \\ [1ex] \hline\hline
  $\lambda$  & Net inflow of susceptibles \\ [1.3ex] \hline
   $\mu$  & Natural death rate \\ [1.3ex] \hline
 $R_0$  & Basic Reproduction Number  \\ [1.3ex] \hline
 $N$  & Total population size \\ [1.3ex] \hline
 $S_0$  & Initial susceptible population   \\ [1.3ex] \hline
\end{tabular}
    \caption{Description of the parameters and quantities related to model \eqref{mod continuo22}. Note that $N=\lambda/\mu$.}
    \label{tab1}
\end{table} 

{\small
\begin{table}[t]
	\centering
 \small
	\begin{tabular}{|l |l |l |}\hline
		Function & Meaning	 &  Assumptions \\ [1ex] \hline\hline
	$A$ & Infectivity function  & $A(t) \geq 0$ \\ [1.3ex] \hline
	$A_{\mu}$ &	\begin{minipage}{5cm} Infectivity function with\\ demography \end{minipage} &   $\hat{A}_{\mu}(0) := \int_0^{\infty}A_{\mu}(t) dt < +\infty$ \\ [1.3ex] \hline
		$\beta$ &	Inhibition function & $\beta(0)=1$, $\beta(M)>0$, $\beta^{\prime}(M) < 0$  \\ [1.3ex] \hline
	$g$ &	Message function    & $g(0)=0$, $g(\cdot) \geq 0$, $g^{\prime}(\cdot) > 0$  \\ [1.3ex] \hline
		$K$ &	Memory kernel    & $K(\tau) > 0$, $\int_0^{\infty} K(\tau) d\tau=1$  \\ [1.3ex] \hline
	%	$S_0$,$F_0$,$M_0$  & Initial values & $S_0 > 0$, $S_0 \leq \lambda /\mu$, $F_0(t) \geq 0$, $M_0(t) \geq 0$\\ [1.3ex] \hline
	\end{tabular}
	\caption{Functions involved in model \eqref{mod continuo22} and the related assumptions.}
	\label{tab:assump}
\end{table}}

We conclude the section with the following:

\begin{prop}\label{prop:1}
	If the IF is exponentially distributed, then model \eqref{mod continuo22} is the behavioural SIR--M model introduced in \cite{donofrio2009}.
\end{prop}

\begin{proof}
	
\noindent Let us take $A(\tau)=\beta_0 e^{-\nu\tau}$, where $\beta_0$ and $\nu$ are two positive constants representing the transmission rate and the average duration of the infectious period, respectively. Set
	
	$$
	I(t) = \int_{0}^{\infty} e^{-(\mu +\nu)\tau} \beta(M(t-\tau)) F(t-\tau) S(t-\tau) d \tau,
	$$

\noindent  substitution of $A_{\mu}(\tau)=\beta_0 e^{-(\mu +\nu) \tau}$ into equation \eqref{BeConstEq} leads to  $F(t) = \beta_0I(t)$. It follows:
	
	\begin{align*}
		\frac{d I(t)}{d t} &= -(\mu +\nu) \int_{0}^{\infty} e^{-(\mu +\nu)\tau} \beta(M(t-\tau)) F(t-\tau) S(t-\tau) d \tau + \beta(M(t)) F(t) S(t) \\
		&= -(\mu +\nu) I(t) + \beta_0 \beta(M(t)) I(t) S(t).
	\end{align*}

\noindent By denoting $\hat{\beta}(M)=\beta_0 \beta(M)$, then

\begin{equation}
	\label{SIRM1}
	\dot{I}= -(\mu +\nu) I + \hat{\beta}(M) S I.
	\end{equation}

\noindent Equations \eqref{BeM},\eqref{EqS} and  \eqref{SIRM1} constitutes the extensively studied behavioural SIR--M model introduced by d'Onofrio and Manfredi  \cite{donofrio2009,donofrio2020,DONOFRIO2022}.
\end{proof}

\begin{remark}
	In the early epidemic model by Kermack and McKendrick \cite{KM1927,KERMACK199133}, it is assumed that factors such as infectivity, recovery, and mortality depend on the age-of-infection, i.e. the time elapsed by the infection. It is worth to mention that the infectious individuals may be further distinguished (say \textit{structured}) according to their age-of-infection. This assumption lead to consider the age-of-infection as an independent variable and the corresponding models are integro-differential models involving partial differential equations. Time-since-infection (and time-since-recovery) structured epidemic models have been studied by several authors (see e.g. \cite{IANNELLIDONOFRIO2023,  Martcheva_2015}). A time-since-infection behavioural model has been recently studied, where the individuals' decisions regarding vaccination is driven by imitation dynamics \cite{Marcheva2016}.    
\end{remark}

\section{Basic properties and existence of equilibria}\label{sec3}
In order to study the solution $\textbf{u}(t)=(S(t),F(t),M(t))$ of model \eqref{mod continuo22}, we assume that, over the interval $(-\infty,0]$, the state variables $S$, $F$ and $M$ are  bounded functions and  $0<S(t)\leq \lambda/\mu$, $F(t)\geq 0$, $M(t)\geq 0$. These properties hold true for all $t \in \mathbb{R}$, in view of the following two theorems.

\begin{theorem}\label{thpos}
    The inequalities $S(t)>0$,  $F(t) \geq 0$, $ M(t) \geq 0$ hold true for all $t > 0$.
\end{theorem}

\begin{proof}
From the first equation of \eqref{mod continuo22} we deduce that

\begin{equation*}
S(t)=e^{-\int_0^t \left[ \mu+\beta(M(\tau))F(\tau)\right]  d \tau} \Big(  S_0 + \int_0^t \lambda e^{\int_0^s \left[ \mu+\beta(M(\tau))F(\tau)) \right] d \tau} ds \Big).
\end{equation*}

\noindent Being $S(0)=S_0 > 0$, it easily follows that $S(t)>0$, $\forall t \geq 0$. The non-negativity of $F(t)$ and $M(t)$ follows from the non-negativity of the functions $A$, $\beta$, $g$ and $K$ (see the assumptions in Table \ref{tab:assump}). 
\end{proof}

\begin{theorem}\label{thbound}
    Assume that 
    $A^{\prime}(\tau) \in L^1([0,\infty))$ and $\lim_{\tau \rightarrow \infty}A(\tau)=0,$ 
    then $S(t) \leq S_{\max}$, $ F(t) \leq F_{\max} $, $ M(t) \leq M_{\max},$ for all $t > 0$, where

\begin{equation}\label{SFMmax}
S_{\max}=\frac{\lambda}{\mu}, \quad
 F_{\max}=\lambda \hat{A}_{\mu}(0)+ \frac{\lambda}{\mu} \int_0^{\infty}  |A_{\mu}'(\tau)| d \tau, \quad
   M_{\max}=g(F_{max}).
\end{equation} 
\end{theorem}

\begin{proof}
The non-negativity of $S(t)$, $F(t)$ and $\beta(M(t))$ implies that $\dot S(t)\leq \lambda-\mu S(t)$. Thus, using the comparison theorem, we obtain

\begin{align*}
 S(t) &\leq e^{-\mu t}\left(S_0 + \lambda \int_0^t e^{\mu s} ds \right)  = \frac{\lambda}{\mu} + \left(S_0-\frac{\lambda}{\mu}\right)e^{-\mu t}.
\end{align*}

\noindent Since $S_0 \leq \lambda/\mu$, we obtain the first bound. From the first and the second equation in \eqref{mod continuo22}, we have

\begin{align*}
    F(t) 
    %&= \int_0^{\infty} \big( \beta(M(t-\tau)) S(t-\tau) F(t-\tau)  \big)  A_{\mu}(\tau)  d \tau \\
   % &= \int_0^{\infty} \big( -\dot S(t-\tau)+\lambda-\mu S(t-\tau) \big)  A_{\mu}(\tau)  d \tau \\
    &= \lambda \int_0^{\infty}  A_{\mu}(\tau)  d \tau - \mu \int_0^{\infty} S(t-\tau) A_{\mu}(\tau) d \tau -\int_0^{\infty} \dot S(t-\tau) A_{\mu}(\tau) d \tau. \\
\end{align*}

\noindent Thus, since $ \int_0^{\infty} S(t-\tau) A_{\mu}(\tau) d \tau \geq 0$,  we deduce that

\begin{equation}\label{F lim}
   F(t)  \leq \lambda \int_0^{\infty}  A_{\mu}(\tau)  d \tau -\int_0^{\infty} \dot S(t-\tau) A_{\mu}(\tau) d \tau.
\end{equation}

\noindent Integrating by parts the second term in the r.h.s of \eqref{F lim}, taking into account the non-negativity of the functions involved and that $\lim_{\tau \rightarrow \infty}A(\tau)=0$, we obtain the following bound for $F(t)$:

\begin{align*}
F(t) &\leq \lambda \hat{A}_{\mu}(0) + \frac{\lambda}{\mu} \int_0^{\infty}  |A_{\mu}'(\tau)| d \tau,
\end{align*}

\noindent where $\hat{A}_{\mu}(0) := \int_0^{\infty}A_{\mu}(\tau) d\tau$. Finally, from the third equation in \eqref{mod continuo22}, since $g$ is increasing and $\int_{0}^{\infty}  K(\tau) d \tau=1$,  we have

\begin{align*}
    M(t) =\int_{0}^{\infty} H(t-\tau) g \big(F(t-\tau)\big) K(\tau) d \tau \leq g(F_{max}) \int_{0}^{\infty}  K(\tau) d \tau = g(F_{max}).
\end{align*}
\end{proof}

%\begin{remark}
    %The condition \eqref{Ainf23} is satisfied if the infectivity function is uniformly continuous.  \textbf{***CHE VUOL DIRE, PERCHE' E PERCHE' INTERESSA?****}
%\end{remark}
	
%\section{Existence of equilibria}\label{sec4}

\noindent The analysis of equilibria begins by observing that model \eqref{mod continuo22} admits a unique \textit{Disease-Free} Equilibrium (DFE) given by 

\begin{equation}\label{dfee}
   DFE = (\Bar{S}_0, 0, 0),
\end{equation}

\noindent where $\Bar{S}_0=\lambda/\mu$. Denote by 

\begin{equation}\label{EEequilibrium}
 EE=(S_{e},F_{e},M_{e}),   
\end{equation} 

\noindent the generic Endemic Equilibrium (EE) of model \eqref{mod continuo22}, i.e. equilibria where $F_{e} > 0$. Endemic equilibria are obtained from conditions

 \begin{equation}\label{eq gen1}
    \begin{array}{ll}
       \lambda-\mu S-\beta(M)SF=0,\quad  \beta(M) S  \hat{A}_{\mu}(0) =1, \quad  
        M = g(F).
    \end{array}
\end{equation}

\noindent Introduce

\begin{equation}\label{R0}
	R_0=\frac{\lambda}{\mu}  \hat{A}_{\mu}(0). 
\end{equation}

%teorema

\begin{theorem}\label{Th stab EE}
  If the assumptions of Theorem \ref{thbound} hold, then: \textit{(i)} if $R_0 <1$, then model \eqref{mod continuo22} has no endemic equilibria;  \textit{(ii)} if $R_0 > 1$, then model \eqref{mod continuo22} has a unique endemic equilibrium.
\end{theorem}

\begin{proof}
   The EE is defined by the following conditions

\begin{equation}\label{SMinfeq}
    M_{e}=g(F_{e}), \quad \quad S_{e}= \frac{\lambda}{\mu R_0 \beta(g(F_{e}))}, 
    \end{equation}

\noindent and $F_{e} > 0$ is the zero of the nonlinear function

    \begin{equation}\label{Finfeq}
    G(x)=x-\mu R_0+\frac{\mu}{\beta(g(x))},
    \end{equation}
    
  \noindent  in the interval $[0,F_{max}]$, where $F_{\max}$ is introduced in \eqref{SFMmax}. 
    We observe that

         \begin{align*}
             &G(0) = -\mu (R_0-1), \\
             &G(F_{max}) = \lambda \int_0^{\infty} |A_{\mu}'(\tau)| d\tau + \frac{\mu}{\beta(g(F_{max}))}>0, \\
             &G'(x)= 1- \frac{\mu \beta'(g(x)) g'(x)}{[\beta(g(x))]^2\hat{A}_{\mu}(0)}>0,
         \end{align*}

\noindent where the last inequality holds since $\beta'(\cdot)<0.$ Thus, when $R_0>1$, it is $G(0)<0$ and $G(x)$ in \eqref{Finfeq} has a unique zero in $[0,F_{\max}].$ Furthermore, if $R_0<1$, then $G(x)$ has no zeros in $[0,F_{max}].$
   
\end{proof}

\begin{remark}
    The quantity denoted by $R_0$ and introduced in \eqref{R0} assumes the meaning of Basic Reproduction Number and it is an important indicator for analyzing the behaviour of epidemic models. In \eqref{R0}, $\hat{A}_{\mu}(0)$ represents the ability to generate secondary cases of infection from an individual within a population (in fact, if $A_{\mu}(\tau)= \beta_0 e^{-(\mu+\nu)\tau}$, the quantity $\hat{A}_{\mu}(0)$ reduces to the transmission rate $\beta_0$ multiplied by the average infectious period $1/(\mu+\nu)$). 

    If $\hat{A}_{\mu}(0)$  is multiplied by the total population size $N=\lambda/\mu$, the meaning of $R_0$, that is the average number of secondary cases produced by one primary infection over the course of the infectious period in a fully susceptible population, is obtained.
\end{remark}

\section{Preliminaries for stability}\label{sec4}

In this section we derive the characteristic equation as preliminary step of the local stability analysis of the equilibria. We use the linearization theory for Volterra integral equations \cite{MILLER1968198}.
We begin by observing that model \eqref{mod continuo22} can be written in integral form as

\begin{align}\label{mod compatto int}
    S(t) &= S_0 + \int_{0}^t \big(\lambda- \mu S(t-\tau) - \beta(M(t-\tau)) S(t-\tau) F(t-\tau) \big) d \tau, \nonumber \\
    F(t) &= F_0(t) + \int_0^t \beta(M(t-\tau)) S(t-\tau) F(t-\tau) A_{\mu}(\tau) d \tau, \\
    M(t) &= M_0(t) + \int_0^t g(F(t-\tau)) K(\tau)  d \tau, \nonumber
\end{align}

\noindent where the functions $F_0(t)$ and $M_0(t)$ are defined by

\begin{equation}
\label{forcingFM}
    \begin{aligned}
        F_0(t) &=\int_{t}^{\infty} \beta(M(t-\tau)) S(t-\tau) F(t-\tau) A_{\mu}(\tau) d \tau,  \\
       M_0(t) &= \int_{t}^{\infty}  H(t-\tau) g(F(t-\tau)) K(\tau)  d\tau = g(F(0)) K(t).
    \end{aligned}
\end{equation}  

\noindent Consider the perturbations  $(W, V, Z)$  to the generic equilibrium $\Tilde{E}=(\Tilde{S}, \Tilde{F}, \Tilde{M})$

\begin{equation*}
S = \Tilde{S} + W, \quad
       F = \Tilde{F} + V, \quad
       M = \Tilde{M} + Z.   
\end{equation*}

\noindent Linearization of \eqref{mod compatto int} about the equilibrium solution leads to the following system:

\begin{equation}\label{mod lineare}
    \textbf{y}(t)=\textbf{f}(t)+ \int_0^t \mathcal{\textbf{Q}} (t-\tau) \textbf{y}(\tau) d\tau,
\end{equation}
where 

\begin{equation*}
    \textbf{y}(t)= \begin{bmatrix}
W(t) \\
V(t) \\
Z(t)
\end{bmatrix}, \quad \textbf{f}(t)= \begin{bmatrix}
W(0) \\
V_0(t) \\
Z_0(t)
\end{bmatrix},
\end{equation*}

\noindent and

\begin{equation*}
  \mathcal{\textbf{Q}} (t)=  \begin{bmatrix}
-\mu -\beta(\Tilde{M}) \Tilde{F} & -\beta(\Tilde{M}) \Tilde{S} & -\beta'(\Tilde{M})\Tilde{F}\Tilde{S}\\ 
\beta(\Tilde{M}) \Tilde{F} A_{\mu}(t) & \beta(\Tilde{M}) \Tilde{S} A_{\mu}(t) & \beta'(\Tilde{M}) \Tilde{S} \Tilde{F} A_{\mu}(t) \\
0 & K(t) g'(\Tilde{F}) & 0
\end{bmatrix}.
\end{equation*}

\noindent Furthermore, the forcing terms are

$$
    V_0(t) = \int_{t}^{\infty} \bigg[ \beta \big(M(t-\tau)\big)  S(t-\tau) F(t-\tau)-\beta(\Tilde{M}) \Tilde{S} \Tilde{F} \bigg] A_{\mu}(\tau) d\tau,
    $$
    $$
    Z_0(t) = \big[g(F(0))-g(\Tilde{F}) \big] K(t). 
$$

\noindent In order to study the stability of $\Tilde{E}$, we refer to an extension of the classical Paley-Wiener theorem. This extension is recalled in  Appendix \ref{AppendixB} (Theorem \ref{pali1}), to which we refer for the meaning of the symbols. System \eqref{mod lineare} fits into the form \eqref{lineare generale}-\eqref{hp lineare gen} with $l=1$, $k=2$, and 

\begin{equation*}
\begin{aligned}
f_0(t) &= W(0), \\
 \textbf{\textit{f}}_1(t) &= {\small\begin{bmatrix} 
 \displaystyle
 %\scriptscriptstyle
 \int_t^{\infty} \big[ \beta \big(M(t-\tau)\big)  S(t-\tau) F(t-\tau) -\beta(\Tilde{M}) \Tilde{S} \Tilde{F} \big] A_{\mu}(\tau) d\tau\\ 
\big[g(F(0))-g(\Tilde{F}) \big] K(t)
\end{bmatrix},} \\
\mathbf{A_{\infty}^1} &= [-\mu -\beta(\Tilde{M}) \Tilde{F}, \quad  -\beta(\Tilde{M}) \Tilde{S}, \quad -\beta'(\Tilde{M})\Tilde{F}\Tilde{S}].
\end{aligned}
\end{equation*}

\noindent Let $\hat{A}_{\mu}$ and $\hat{K}$ denote the Laplace transforms

\begin{equation}\label{AK hat}
\hat{A}_{\mu}(w)=\int_0^{\infty}e^{-w t} A_{\mu}(t)dt , \quad \hat{K}(w)=\int_0^{\infty}e^{-w t}K(t)dt,
\end{equation}

\noindent with $w \in \mathbb{C}$. If $Re(w)>0$, the characteristic equation of system \eqref{mod lineare}  takes the form

\begin{equation}\label{det ric}
 det\begin{bmatrix}
    w+\mu +\beta(\Tilde{M}) \Tilde{F} & \beta(\Tilde{M}) \Tilde{S} & \beta'(\Tilde{M})\Tilde{F}\Tilde{S} \\ 
-\beta(\Tilde{M}) \Tilde{F}\hat{A}_{\mu}(w) & 1-\beta(\Tilde{M}) \Tilde{S} \hat{A}_{\mu}(w) & -\beta'(\Tilde{M}) \Tilde{S} \Tilde{F} \hat{A}_{\mu}(w) \\
0 & -g'(\Tilde{F}) \hat{K}(w)  & 1
\end{bmatrix}=0.
\end{equation}

\noindent Under the assumptions of Theorem \ref{thbound}, we have that the functions $f_0$ and $ \textbf{\textit{f}}_1$ satisfy the hypothesis of Theorem \ref{pali1}. The following theorem gives a sufficient condition for the local asymptotic stability of the equilibria.

\begin{theorem}\label{th Stab cont} Any equilibrium $\Tilde{E}=(\Tilde{S},\Tilde{F}, \Tilde{M})$ of model \eqref{mod continuo22} is locally asymptotically stable if the characteristic equation 
	\begin{equation}\label{Eq car th}
		 \beta(\Tilde{M}) \Tilde{S} \hat{A}_{\mu}(w) \bigg( 1+ \frac{g'(\Tilde{F})\beta'(\Tilde{M})\Tilde{F} \hat{K}(w)}{\beta(\Tilde{M})} \bigg) = 1+ \frac{\beta(\Tilde{M}) \Tilde{F}}{w + \mu},
	\end{equation}
 
	\noindent has no solution for $ Re(w) \geq 0$.
 
\end{theorem}
\begin{proof} \textit{(Sketch)}
    The proof is a straightforward application of Theorem \ref{pali1}, after observing that the transcendental equation \eqref{Eq car th} is equivalent to \eqref{det ric}. 
\end{proof}

\noindent We conclude the section by noting that  the information dependent parameters affect the location of the  equilibrium EE only through the message function (see the term 
$g(F_{e})$ in \eqref{SMinfeq}, \eqref{Finfeq}), while its stability properties depend also on the memory kernel $K$. 

%\begin{remark}
   % In the next, we focus only on the scenario where $Re(w) > 0$. However, when $Re(w)=0$, it's easy to adapt the following proof, by interpreting the equation as the limit for $\Tilde{w} \rightarrow w$, with $Re(\Tilde{w}) >0$.  
%\end{remark}
	
\section{Local stability analysis of equilibria}\label{sec5}

\noindent In this section we analyze the stability properties of the equilibria DFE and EE, by applying Theorem \ref{th Stab cont} to some special cases of interest. In our analysis we assume that $Re(w)>0$. All the results can be easily extended to $Re(w) \geq 0$, by passing to the limit for $\Tilde{w} \rightarrow w$ with $Re(\Tilde{w})>0$. 

At the DFE \eqref{dfee} the characteristic equation \eqref{Eq car th} simply reduces to $\Bar{S}_0 \hat{A}_{\mu}(w) = 1$ and  $Re(w) > 0$ implies

$$
\bigg| \frac{\lambda}{\mu} \hat{A}_{\mu}(w) \bigg| \leq \frac{\lambda}{\mu} \hat{A}_{\mu}(0).
$$

\noindent Therefore, from Theorem \ref{th Stab cont} and \eqref{R0} the DFE is locally asymptotically stable if $R_0 < 1$.

Regarding the equilibrium EE \eqref{EEequilibrium}, from the characteristic equation \eqref{Eq car th} we obtain sufficient stability conditions in some specific cases. Generally speaking, the study of the roots of the characteristic equation is based on the following considerations: when $Re(w)>0$ it is 

\begin{equation}\label{condition 1eq car}
     |\beta(M_{e}) S_{e} \hat{A}_{\mu}(w)| \leq 1, \quad  Re \bigg( 1+ \frac{\beta(M_{e}) F_{e}}{w + \mu} \bigg) >1,
\end{equation}
where the first inequality comes from condition \eqref{eq gen1}$_2$  and the second one is a direct consequence of $Re(w) > 0$.
 Given \eqref{condition 1eq car}, a sufficient condition guaranteeing that the  characteristic equation \eqref{Eq car th} does not have roots $w$ such that $Re(w) \geq 0$ (as required by Theorem \ref{th Stab cont}) is given by  

 \begin{equation}\label{suff cond }
      \bigg| 1+ \frac{g'(\Tilde{F})\beta'(\Tilde{M})\Tilde{F} \hat{K}(w)}{\beta(\Tilde{M})} \bigg|\leq 1.
 \end{equation}

\noindent  In the following, we discuss three noteworthy cases according to the form of the memory kernel $K$: no-delay; exponentially fading memory; unimodal memory.

\medskip

\noindent \textit{Case (i): no-delay (Dirac's delta memory kernel, $K(\tau)=\delta(\tau)$)}. From \eqref{Eq car th} we get

    \begin{equation}\label{eq car31} g'(F_{e})\beta'(M_{e})F_{e} S_{e} \hat{A}_{\mu}(w)  + \beta(M_{e}) S_{e} \hat{A}_{\mu}(w) = 1 + \frac{\beta(M_{e}) F_{e}}{w + \mu}.
    \end{equation}

    \begin{theorem}
    \label{t1}
           If $K(\tau)=\delta(\tau)$ and $A_{\mu}(\tau)$ is a positive definite function\footnote{A continuous function $a: \mathbb{R}\rightarrow \mathbb{C}$ is called a positive definite function if $\sum_{i,j}a(x_i-x_j)z_i \bar{z}_j \geq 0$ for any choice of finite sequences $\{x_n\}_n$, $\{z_n\}_n$, with $x_n \in \mathbb{R}$ and $z_n \in \mathbb{C}$ (see \cite{LUBICH1983}).}, then EE is LAS for any choice of the functions $\beta$ and $g$.
        \end{theorem}
        \begin{proof} 
            Let be $Re(w) > 0$. For the Bochner characterization of positive definite functions (see Section 8.2 in \cite{LUBICH1983}), it is
            $Re \big( \hat{A}_{\mu}(w) \big) \geq 0.$
             
          Thus, in view of \eqref{condition 1eq car}, since $g'(\cdot)>0$ and $\beta'(\cdot)<0$, the real part of the left-hand side of \eqref{eq car31} is less than 1, while the real part of the r.h.s. is strictly greater than 1. Therefore  equation \eqref{eq car31} can never be satisfied.
        \end{proof}
    
    \begin{remark}
        An example of positive definite infectivity function is $A_{\mu}(\tau)=e^{-(\mu + \nu)\tau}$, for which model \eqref{mod continuo22} is the $SIR-M$ model \eqref{EqS}-\eqref{SIRM1}, as shown in Proposition \ref{prop:1}. More generally, any non-negative, non-increasing infectivtiy function, convex on $[0,\infty)$, is positive definite \cite{LUBICH1983} and thus Theorem \ref{t1} applies.
        \end{remark}
   
   \smallskip
    
   \begin{theorem}
     If $K(\tau)=\delta(\tau)$,  $A'(\tau)<0$, and the functions $\beta$ and $g$ are such that
     
        $$g'(F_{e})\beta'(M_{e})F_{e} \leq - \beta(M_{e}),$$
         then EE is LAS.
   \end{theorem}
   \begin{proof}
       Let be $\omega=x+iy$, with $x>0.$ We treat the three cases $y=0,$ $y>0,$ and $y<0$ separately.
\begin{itemize}
\item If $y=0,$ we consider equation \eqref{eq car31} for which the inequalities \eqref{condition 1eq car} hold. Since in this case $g'(F_{e})\beta'(M_{e})F_{e} S_{e} \hat{A}_{\mu}(w)<0,$ for the assumptions on $\beta$ and $g$ in Table \ref{tab:assump}, the left-hand-side is less than one and cannot be equal to the r.h.s., which is in turn greater than one. So, equation \eqref{eq car31} has no roots with positive real part. 
\item If $y>0,$  we write the imaginary part of equation \eqref{eq car31} as

  \begin{equation}\label{car immag}
           \bigg( 1+ \frac{g'(F_{e})\beta'(M_{e})F_{e}}{\beta(M_{e})} \bigg) y \int_0^{\infty} e^{-x t} \sin(yt) A_{\mu}(t) dt = \frac{\beta(M_{e}) F_{e} y^2}{(x + \mu)^2+y^2}. 
       \end{equation}

Inspired by the proof of 
\cite[Lemma 1]{Driessche_Watmough_2000} we 
observe that,  since, for the theorem 
assumption, $e^{-xt}A_{\mu}(t)$ is a 
strictly decreasing function, it is 

{\small
 $$
 \int_0^{\infty}e^{-xt}\sin{(yt)}A_{\mu}(t)dt=  \sum_{j=0}^{+\infty}\int_{\frac{2\pi j}{y}}^{\frac{2\pi (j+1/2)}{y}}(e^{-xt}A_{\mu}(t)-e^{-x(t+\pi/y)}A_{\mu}(t+\pi/y))\sin{(yt)}dt>0.$$}

It follows that 

\[y \int_0^{\infty} e^{-x t} \sin(yt) A_{\mu}(t) dt >0.\]

The positivity of the r.h.s. of \eqref{car immag} implies that  the 
characteristic equation \eqref{eq 
car31} can never be satisfied if 

$$
       1+ \frac{g'(F_{e})\beta'(M_{e})F_{e}}{\beta(M_{e})}\leq 0.
       $$
 
\item 
If $y=-z<0,$ then $\sin{(yt)}=-\sin{(zt)}$. This case is treated as the previous one since

$$\int_0^{\infty}e^{-xt}\sin{(yt)}A_{\mu}(t)dt=-\int_0^{\infty}e^{-xt}\sin{(zt)}A_{\mu}(t)dt<0.$$
\end{itemize}
   \end{proof}

\medskip

\noindent \textit{Case (ii): exponentially fading memory (weak Erlang memory kernel, $K(\tau)= a e^{-a \tau}$).} The Laplace transform in \eqref{AK hat} reads $\hat{K}(w)=a/(a+w)$. Therefore,  
 the sufficient condition \eqref{suff cond } can be written as
 
 \begin{equation}\label{suff cond weak}
    \bigg| 1+ \frac{a g'(F_{e})\beta'(M_{e})F_{e}}{\beta(M_{e})(a+w)} \bigg| \leq 1. 
\end{equation}

 \begin{theorem}
 \label{t4}
   If $K(\tau)= a e^{-a \tau}$ and the functions $\beta$ and $g$ satisfy the inequality
    
  \begin{equation}\label{hypTh5.4}  g'(F_{e})\beta'(M_{e})F_{e} \geq -2 \beta(M_{e}),
\end{equation}
    then EE is LAS for any choice of $A_{\mu}(\tau)$.
\end{theorem}
\begin{proof}
   Let us set $w=x+iy$,  with $x>0$, and $c=\frac{g'(F_{e})\beta'(M_{e})F_{e}}{\beta(M_{e})}$. It is

\begin{equation*}
   \bigg| 1+ \frac{a g'(F_{e})\beta'(M_{e})F_{e}}{\beta(M_{e})(a+w)} \bigg|^2 = 1+ \frac{c a }{(a+x)^2+y^2} \bigg( ca + 2(a+x) \bigg). 
\end{equation*}

\noindent Therefore, the sufficient condition \eqref{suff cond weak} is satisfied if and only if

$$
-1 \leq \frac{c a }{(a+x)^2+y^2} \bigg( ca + 2(a+x) \bigg) \leq 0.
$$
Taking into account that $\beta'(\cdot)<0$, this condition is always satisfied if $- c \leq 2$, which is equivalent to \eqref{hypTh5.4}.

%holds true. In this case the left-hand-side of \eqref{p3} has modulus less than or equal to one (see the first of \eqref{condition 1eq car}) and the real part of the r.h.s. is greater than one, thus excluding the existence of roots for $Re(w) > 0$.
\end{proof}

\begin{remark}
    The result given in 
    Theorem \ref{t4} remains valid in the no-delay case  $K(\tau)=\delta(\tau)$.
    As matter of fact, the sufficient condition \eqref{suff cond } turns out to be

    \begin{equation*}
    \bigg| 1+ \frac{g'(F_{e})\beta'(M_{e})F_{e} }{\beta(M_{e})} \bigg| \leq 1,
    \end{equation*}
    which is equivalent to \eqref{hypTh5.4}, since $\beta'(M_{e})<0$.
\end{remark}

\medskip

\noindent \textit{Case (iii): unimodal memory (strong Erlang memory kernel, $K(\tau)= a^2 \tau e^{-a \tau}$).} The Laplace transform in \eqref{AK hat} reads $\hat{K}(w)=  \frac{a^2}{(a+w)^2} $. Therefore,  the sufficient condition \eqref{suff cond } is
 
 \begin{equation}\label{strongdim1}
    \bigg| 1+ \frac{a^2 g'(F_{e})\beta'(M_{e})F_{e}}{\beta(M_{e})(a+w)^2} \bigg| \leq 1.
    \end{equation}

\noindent Reasoning as in the proof of Theorem \ref{t4}, let us set $w=x+iy$, with $x>0$, and $c=\frac{g'(F_{e})\beta'(M_{e})F_{e}}{\beta(M_{e})}$. It is

\begin{align*}
   \bigg| 1+ \frac{\hat{K}(w) g'(F_{e})\beta'(M_{e})F_{e}}{\beta(M_{e})} \bigg|^2 &= 1 + 2 c a^2 \frac{a^2((a+x)^2-y^2)}{((a+x)^2+y^2)^2}+ \frac{c^2a^4}{((x+a)^2+y^2)^2}.\\ 
\end{align*}
Since $\beta'(\cdot)<0$, it follows that \eqref{strongdim1} is satisfied if and only if 

$$
-1 \leq  2 c a^4 \frac{((a+x)^2-y^2)}{((a+x)^2+y^2)^2}+ \frac{c^2a^4}{((x+a)^2+y^2)^2} \leq 0.
$$

\noindent We observe that, while the left inequality is obviously accomplished, the right one can be written as

$$
2(a+x)^2 \geq -c a^2+ 2y^2,
$$

\noindent which is not true for   any $w$ with $Re(w) \geq 0$, as required by Theorem \ref{th Stab cont}. Therefore, instability of the equilibrium EE may take place in principle. In the next section, through numerical simulation, we will show that self-sustained oscillations may occur indeed. 

\section{Applications}\label{sec6}
%We apply the theoretical results obtained in the previous sections to the case of specific infectious diseases for which indications about the functions and values of the parameters to be used are available in the literature.

\subsection{Preliminaries}

The integro-differential system \eqref{mod continuo22} cannot be completely solved by analytical methods and thus numerical simulations are essential for gaining insights into the qualitative properties of the solutions and for providing quantitative assessment on relevant quantities, like, e.g., epidemic peaks. In the literature, this problem has been faced for renewal equations describing epidemics where the demographic turnover and the behavioural response of the population are neglected
 \cite{MPV2022,ELE2023}. In such cases a non-standard method was built \textit{ad hoc} for the problem at hand, in order to preserve the positivity of the solution and the asymptotic behaviour of the continuous model, without any restriction on the discretization stepsize.
 Here, we build a new version of the non-standard method to deal with the behavioural integral epidemic model with demographic turnover \eqref{mod continuo22}. Some details on the discretization scheme are reported in Appendix \ref{AppendixA}.
 
 %For this purpose, for any fixed value of the stepsize $h>0,$ we have treated (viewed) the numerical method as a discrete-time RE model and  obtained results on the stability of numerical equilibria which replicate the ones described in Section 6 for the continuous model. Having shown the existence of numerical equilibria and  studied their stability, we  have studied their approximation properties proving the convergence as $h\to 0.$ The resulting numerical method performs efficient and reliable  long-time simulations, which are essential to gain a deep understanding in equilibria of the system and their stability properties related to control parameters.

We finally remark that model \eqref{mod continuo22} requires  the initial condition $S_0=S(0)$ and the forcing functions $F_0(t)$ and $M_0(t)$ defined in \eqref{forcingFM}, which depend on the entire past history of $S,F$ and $M$ over $(-\infty,0]$. Here, we assume that the reaction of the population to the epidemic is triggered at time $t=0$. Therefore, the initial condition $S_0$ is the steady state of the susceptible population for ``non-- behavioural''  case $M=0$.

The constitutive law \eqref{BeConstEq} may be written as

\begin{equation*}
    F(t) =\int_{0}^{\infty} \beta(M(t-\tau))  S(t-\tau) F(t-\tau)  H(t-\tau) A_{\mu}(\tau) d \tau, 
    \end{equation*}
 where $H(\cdot)$ is the Heaviside step function. Assuming $\dot S(\tau)\approx 0$, for $\tau\leq 0,$ one gets

\begin{equation}\label{ForcingF00}
    F_0(t) =  A_{\mu}(t)(\lambda-\mu S_0). 
    \end{equation}

\noindent In the next subsections we will numerically test model \eqref{mod continuo22} by choosing some different infectivity functions taken from studies where specific diseases like Influenza and SARS are considered \cite{Aldis2005,donofrio2009,glass2009,Roberts2006}.

\subsection{Exponential and unimodal IF}

We begin with the exponentially distributed IFD already considered in Proposition \ref{prop:1}:

\begin{equation}\label{IFD1}
	A_{\mu}(\tau)=\beta_0 e^{-(\mu+\nu) \tau}.  
\end{equation}

\noindent As shown in Proposition \ref{prop:1}, in this case model \eqref{mod continuo22} becomes the SIR--M model studied in 2009 by d'Onofrio and Manfredi to investigate steady oscillations exhibited by a variety of diseases like pre-vaccination measles \cite{donofrio2009}. Here, and in all the subsequent subsections, we  choose the message function and the inhibition function as follows: 

\begin{equation*}
g(F)=F, \quad \quad \beta(M)=\frac{1}{1+\alpha M},
\end{equation*}

\noindent where $\alpha>0$ is a decline factor due to voluntary social distancing.

We numerically test the model by using the parameter constellation \cite{donofrio2009}:  $R_0=20$, $\nu=1/7\text{ days}^{-1}$, $\mu=1/75\text{ years}^{-1}$, $\alpha=8\cdot10^3$, and $N=5\cdot10^7$, where the values of $R_0$ and $\nu$ refer to measles disease. \\
 It is well known that the endemic equilibrium of SIR--M model is locally asymptotically stable when the memory kernel is a weak Erlang kernel  \cite{donofrio2009,wang2006}. 
 This kernel satisfies the sufficient condition for the 
 asymptotic local stability given in Theorem \ref{t1} and \ref{t4}.

\begin{figure}
	\centering
	\begin{subfigure}[b]{0.47\textwidth}
		\centering
		\includegraphics[width=\textwidth]{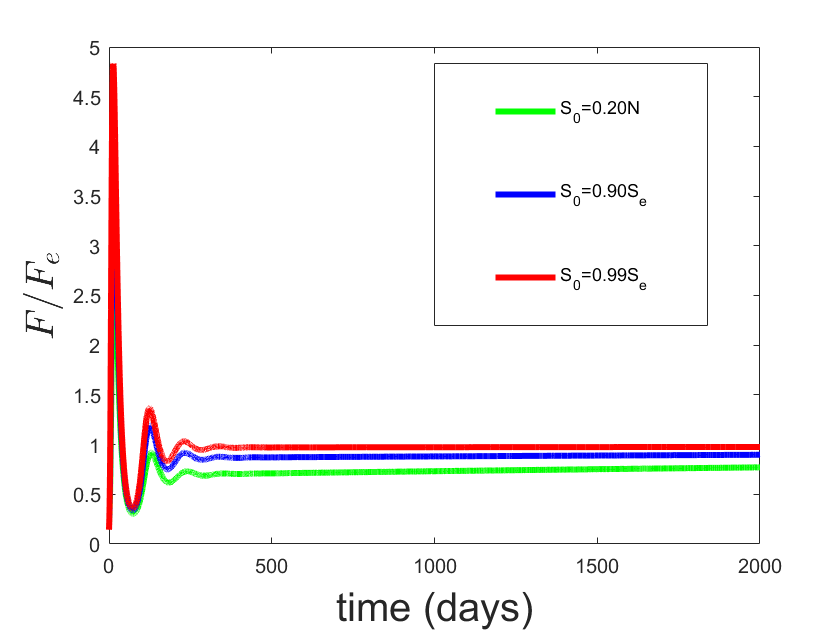}
	\end{subfigure}
	\hfill
	\begin{subfigure}[b]{0.47\textwidth}
		\centering
		\includegraphics[width=\textwidth]{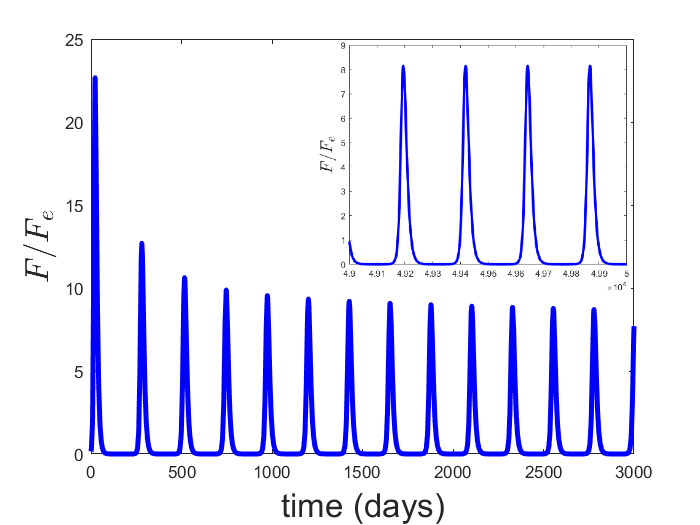}
	\end{subfigure}
	\hfill
			\caption{Time profile of the ratio $F/F_{e}$ as predicted by model \eqref{mod continuo22} with exponential IFD \eqref{IFD1}. Left panel: case with  with weak Erlang memory, i.e. $K(\tau)=a e^{-a \tau}$.    Three different initial conditions are reported:  $S_{0}=0.20N$, $S_{0}=0.90S_{e}$ and $S_{0}=0.99S_{e}$, where  $S_e=0.2648.$
           Right panel: case with strong Erlang memory kernel, i.e. $K(\tau)=a^2 \tau e^{-a \tau}$  (initial condition  $S_0= 0.2N$). Parameter values: see the text.}
			\label{fig:1}
\end{figure}

The numerical simulations enable the visualization of the theoretical results. Taking  $K(\tau)=a e^{-a \tau}$ with $a=1/30\text{ days}^{-1}$,  we observe in Figure \ref{fig:1}, left panel, an asymptotic convergence to the endemic equilibrium. In order to show the algorithm's robustness with respect to perturbations of the initial value, we consider three different scenarios:  a generic initial value, say $S_{0}=0.20N$; high and low perturbation of the EE, i.e. $S_{0}=0.90S_{e}$ and $S_{0}=0.99S_{e}.$ In all of these scenarios, the stability result is confirmed, although perturbations far from the endemic equilibrium show a  slow convergence (see the green line in Figure \ref{fig:1}, left panel).

\begin{figure}
	\centering
	\begin{subfigure}[b]{0.47\textwidth}
		\centering
		\includegraphics[width=\textwidth]{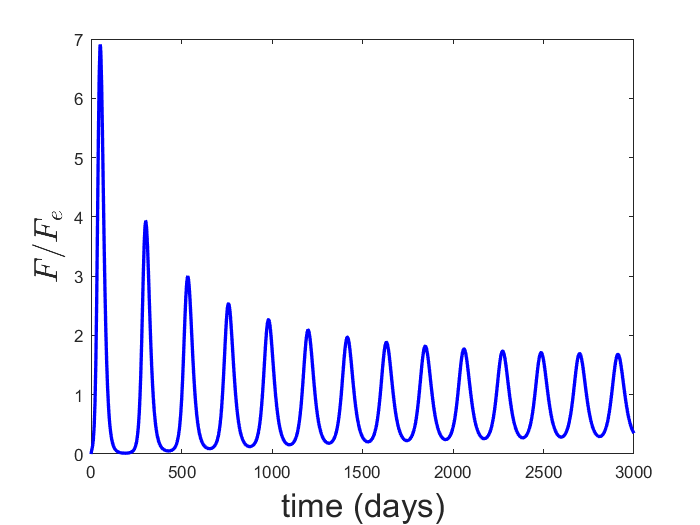}
	\end{subfigure}
	\hfill
	\begin{subfigure}[b]{0.47\textwidth}
		\centering
		\includegraphics[width=\textwidth]{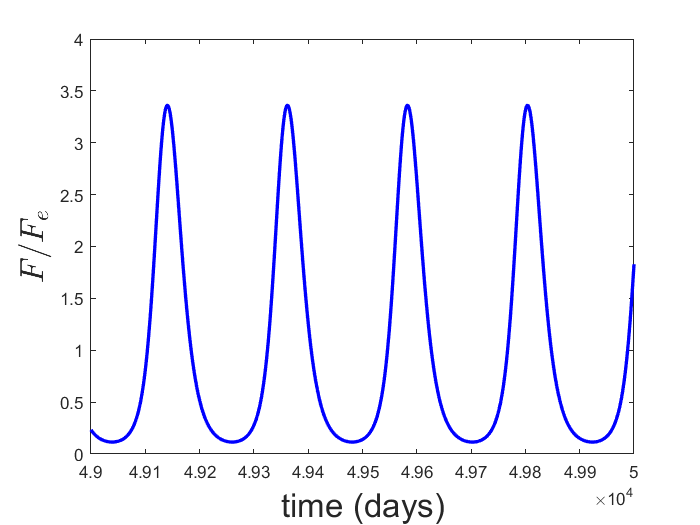}
	\end{subfigure}
	\hfill
			\caption{Time profile of the ratio $F/F_{e}$ as predicted by model \eqref{mod continuo22}.  Case with exponential IFD and strong Erlang memory kernel, i.e. $K(\tau)=a^2 \tau e^{-a \tau}$  \eqref{IFD3} (initial condition  $S_0= 0.2N$).  Parameter values: see the text. }
			\label{fig:2}
\end{figure}

 Now, let us assume that the memory kernel has a strong Erlang distribution $(n=2)$ $K(\tau)=a^2 \tau e^{-a \tau}.$
 Under such hypothesis it has been shown that, when the IFD is given by \eqref{IFD1}, the endemic equilibrium of the associated SIR--M model can be  destabilized via Hopf bifurcations yielding to stable recurrent oscillations \cite{donofrio2009,DONOFRIO2022}. We recall that the sufficient conditions \eqref{strongdim1} guaranteeing the local asymptotic stability of the endemic equilibrium cannot be satisfied in this case. Indeed, oscillations may occur (see Figure \ref{fig:1}, right panel).
In order to analyse the nature of these oscillations, we consider, in a given interval $[T_1,T_2]$, the amplitudes $A_1$ and $A_2$ of the first and last oscillation, respectively. We introduce the percentage \emph{amplitude relative variation} $A_{var}$ as 
 \[A_{var}=\frac{A_2-A_1}{A_2}\cdot 100.\]
 We get $A_{var}\approx 0.018\%$  in the time interval $[4.9\cdot 10^4, 5.0\cdot 10^4]$ (see 
 Figure \ref{fig:1}, right panel, top right square).  Therefore, barring numerical errors, oscillations appear to be self-sustained.
 The integral formulation of the renewal equation \eqref{mod continuo22} allows to verify that the occurrence of oscillations can be observed not only when the IFD is given by \eqref{IFD1},  but also in unimodal cases.  For example a similar result may be observed when the   IFD is given by:

\begin{equation}\label{IFD3}
	A_{\mu}(\tau)=\beta_0  \tau e^{-(\mu+\nu) \tau}.
\end{equation}

\noindent In this case, we still observe oscillations with 
$A_{var}\approx 0.017\%$ in the time interval  $[4.9\cdot 10^4, 5.0\cdot 10^4]$ (see in  Figure \ref{fig:2}).

\subsection{Trapezoidal IF}\label{SARS sub}

\noindent In 2006, M. G. Roberts proposed an integral equation model for describing the transmission of a generic infectious disease invading a susceptible population \cite{Roberts2006}. The model was motivated by the global epidemic risk of SARS in 2003. Roberts assumed that the IF, $A(\tau),$ can be represented by the following trapezoidal distribution, requiring minimal data for its description:

\begin{equation}\label{infect sars}
    A(\tau)= \begin{cases}
         p_0 \frac{\tau-\tau_a}{(\tau_b-\tau_a)} & \tau_a < \tau < \tau_b \\
        p_0 & \tau_b< \tau < \tau_c \\
         p_0\frac{\tau_d-\tau}{(\tau_d-\tau_c)} & \tau_c < \tau < \tau_d \\
         0 & \text{otherwise}
     \end{cases}.
 \end{equation}
In \eqref{infect sars}, the parameter $p_0$ weights the contacts between the infected and susceptible individuals. It is assumed that all infected individuals have the ability to transmit the infection between time intervals $\tau_b$ and $\tau_c$ after their own infection, while no transmission occurs before time $\tau_a$ and after time $\tau_d$. The transmission probabilities between time intervals $\tau_a$ and $\tau_b$, as well as between times $\tau_c$ and $\tau_d$, are estimated through linear interpolation.  Function \eqref{infect sars} was introduced for the first time by Aldis and Roberts in order to study potential smallpox epidemics \cite{Aldis2005}. 

Here we employ the IF \eqref{infect sars} into the integral model \eqref{mod continuo22}. We also assume that a strong Erlang distribution describes the  memory of the population, $K(\tau)=a^2 \tau e^{-a \tau}.$
  Such a kernel means that the current information is un-available and the maximum weight is assigned at the information arrived to the public after a characteristic time $T = 2/a$. This appears to be somehow coherent with the fact that the available information about SARS was scarce at the time it broke out. 
 
  According to \cite{Roberts2006}, we take: $R_0=3.3$ and $(\tau_a,\tau_b,\tau_c,\tau_d)=(4,7,11,14)$ days. We also take the information parameter values $\alpha=5 \cdot 10^4$, $a=1/30$ and the demographic parameter values $S_0=0.7N$, $\mu=1/75\text{ years}^{-1},$ and $N=5\cdot10^7$. The parameter $p_0$, which can be obtained from \eqref{R0}, is $p_0=9.4 \cdot 10^{-9}$.
   
\begin{figure}
	\centering
	\begin{subfigure}[b]{0.47\textwidth}
		\centering
		\includegraphics[width=\textwidth]{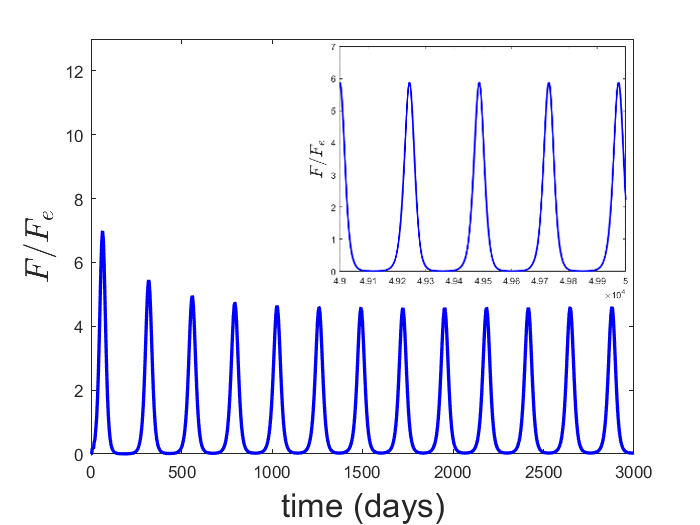}
	\end{subfigure}
	\hfill
	\begin{subfigure}[b]{0.47\textwidth}
		\centering
		\includegraphics[width=\textwidth]{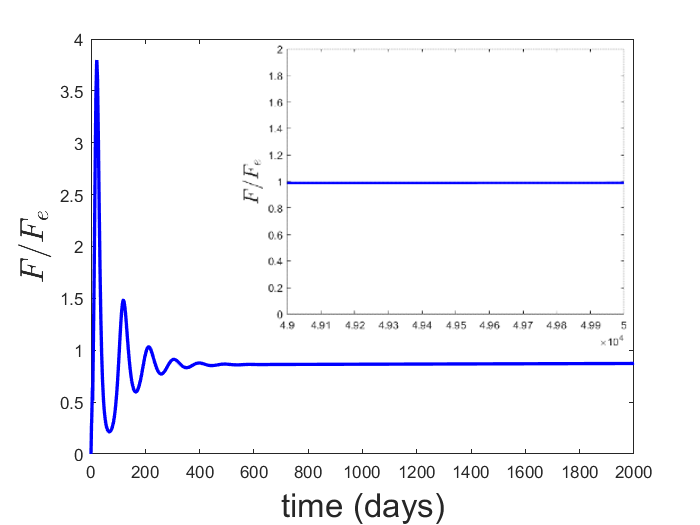}
	\end{subfigure}
	\hfill
	\caption{Time profile of the ratio $F/F_{e}$ as predicted by model \eqref{mod continuo22}.  Left panel: case with strong Erlang memory kernel, i.e. $K(\tau)=a^2 \tau e^{-a \tau}$ and IF given by \eqref{infect sars}.  Right panel:  case with weak Erlang memory kernel, i.e. $K(\tau)=a e^{-a \tau}$ and IF given by \eqref{infect influ}. The initial condition is $S_0=0.7N$. Parameter values: see the text.  }
	\label{fig:3}
\end{figure}

 In Figure \ref{fig:3}, left panel, it is shown that the interplay between trapezoidal distribution \eqref{infect sars} and the strong Erlang memory kernel produces  oscillations that appear to be self-sustained since in the time interval $[4.9\cdot 10^4,5.0\cdot 10^4]$ we get $A_{var}\approx -0.009\%$ (see Figure \ref{fig:3}, left panel, top right square).

%\subsection{Influenza}\label{influ}
In 2007, Roberts and coworkers \cite{Roberts2007} studied the transmission of influenza within a population by using the following IF:

\begin{equation}\label{infect influ}
    \tilde{A}(\tau)= 
         \frac{1}{T_I} A(\tau), 
 \end{equation}

\noindent where $A(\tau)$ is given by \eqref{infect sars} and $T_I$ is the mean
infectious period, $T_I=(\tau_d+\tau_c-\tau_a-\tau_b)/2$. 

%We can also make the assumption that $T_I$ corresponds to the time when symptoms first appear, as proposed by Glass and Becker in \cite{glass2009}. 

Here, we assume that a weak Erlang distribution describes the  memory of the population, $K(\tau)=a e^{-a \tau}.$
  Such a kernel means that  the maximum weight is given to the current information. 
  
  According to \cite{Roberts2007}, we take: $R_0=3$ and $(\tau_a,\tau_b,\tau_c,\tau_d)=(1.2, 2.0, 5.3, 6.1)$ days. All the other parameter values are as in the previous experiment. The parameter $p_0$, which can be obtained from \eqref{R0}, is $p_0= 6.0 \cdot 10^{-8}$.

  %Influenza p0=6.001391036040249e-08
   
 In Figure \ref{fig:3}, right panel, it can be observed that the
  fraction $F/F_{e}$ slowly converges after dumped oscillations, in agreement with Theorem \ref{t4}. 

\subsection{A birth-death process}
Glass and Becker studied the influenza transmission within households by assuming that changes in the infectiousness of individuals may be described by a deterministic birth–death process for the growth of the virus population \cite{glass2009}. The death rate applies once the immune system becomes active, $T_I$
days after infection, where $T_I$ is assumed to correspond to the time when symptoms first appear. The proposed infectiousness profile can be obtained, for example, by the following IF:

\begin{equation}\label{influ house}
     A(\tau)= \begin{cases}
         p_0 e^{-2}e^{2.5\tau} &  \tau \leq T_I \\
        p_0 e^{5}e^{-\tau} &  \tau > T_I
     \end{cases},
\end{equation}

\begin{figure}
	\centering
	\begin{subfigure}[b]{0.47\textwidth}
		\centering
		\includegraphics[width=\textwidth]{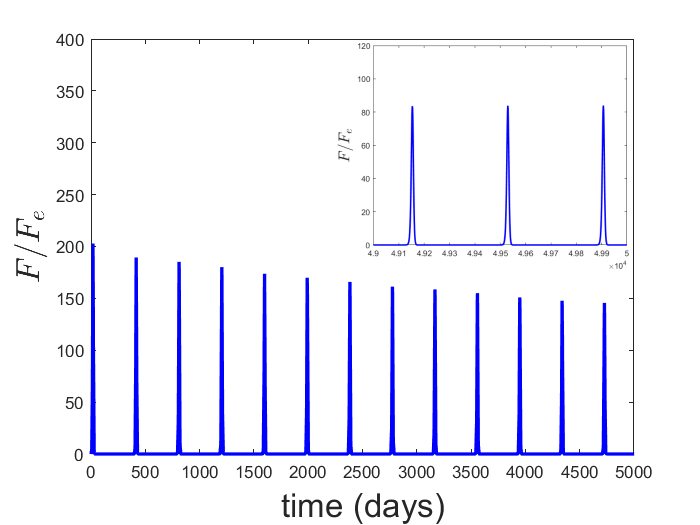}
	\end{subfigure}
	\hfill
	\begin{subfigure}[b]{0.47\textwidth}
		\centering
		\includegraphics[width=\textwidth]{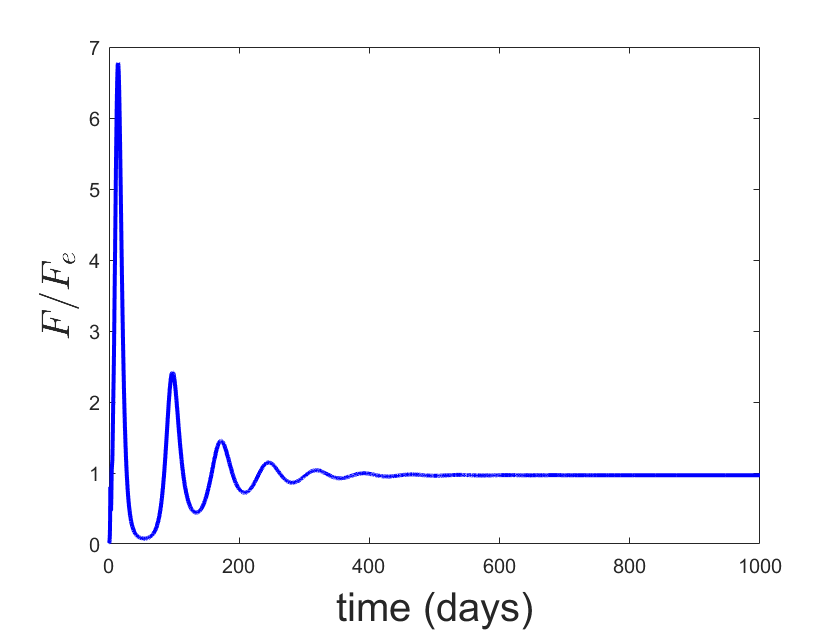}
	\end{subfigure}
	\hfill
	\caption{Time profile of the Force of infection $F/F_e$ according to model \eqref{mod continuo22}. Left panel:  case with strong Erlang memory kernel, i.e. $K(\tau)=a^2 \tau e^{-a \tau}$ and IF given by \eqref{influ house}. Right panel:  case with weak Erlang memory kernel, i. e. $K(\tau)=a e^{-a \tau}$ and IF given by \eqref{influ house}. The initial condition is $S_0=0.98N$. Parameter values: see the text.}
	\label{fig:4}
\end{figure}

\noindent We take $T_I=2$ and $R_0=3$.  We also take the information parameter values $\alpha=8 \cdot 10^4$, $a=1/30$ and the demographic parameter values $S_0=0.98S_e$, $\mu=1/75\text{ years}^{-1},$ and $N=5\cdot10^7$.  The parameter $p_0$, which can be obtained from \eqref{R0}, is $p_0= 2.14 \cdot 10^{-9}$. In Figure \ref{fig:4}, it is shown the time profile of the quantity $F/F_{e}$: in the left panel we can see the epidemic waves produced by the strong Erlang memory, $K(\tau)=a^2 \tau e^{-a \tau}$. We observe oscillations with very large  amplitude.  In the time interval $[4.9\cdot 10^4, 5.0\cdot 10^4]$ it is $A_{var}= 0.5\%$
(see Figure \ref{fig:4}, left panel,  top right square).

\noindent In the right panel dumped oscillations towards a stable equilibrium may be seen in case of weak Erlang memory, $K(\tau)=a e^{-a \tau}$.

%il p0 per l'influenza negli household è 2.200377645738645e-09

\section{Conclusions and future perspectives}\label{sec7}

In the field of ME, a vast portion of literature is based on the use of compartmental epidemic models that remain a  natural extension of the SIR compartmental model proposed by Kermack and McKendrick. However, the early model credited to Kermack and McKendrick is well recognized as a {\em renewal} equation with an integral formulation that is more general than many compartmental models ruled by ODEs. Consequently, before building a compartmental model, it may be advisable to first attempt to formulate the problem in terms of an integral renewal equation \cite{DIEKIN}.  This formulation is very significant for models incorporating the opinion-driven change of human behaviour. In fact, behavioural effects can often serve as the primary determinant in shaping the progression of an epidemic. In particular, how the population keeps memory of information and rumours regarding the disease may influence the solution of the epidemic models, thereby inducing, for example, recurrent oscillations which, in turn, may be interpreted as multiple epidemic waves. 

In this paper, we have proposed an integral renewal equation describing the epidemic propagation of infectious diseases. The equation is {\it behavioural} in the sense that the force of infection includes the information index to account for the “human element”, namely how individuals modify their risky contact according to information and rumours regarding the disease.
The main results are the following:

(\textit{a}) we obtain sufficient conditions, written in terms of the information parameters and functions, for the endemic equilibrium of the behavioural integral epidemic model to be locally stable. In particular, we show that the model admits a stable endemic equilibrium when the memory kernel is a weak Erlang kernel; 

(\textit{b}) with the aid of numerical simulations based on a non-standard numerical method, we show that the model may admit self--sustained oscillations when the memory is more focused in the disease's past history, as exemplified by the strong Erlang kernel. This result is in agreement with the analysis provided by A. d'Onofrio and P. Manfredi, who proved that oscillations may be induced by information-related changes in contact patterns when the disease spread is described by a SIR model with information--dependent contact pattern \cite{donofrio2009};

(\textit{c}) we also consider some different infectivity functions used in the literature to describe infectious disease like Influenza and SARS. Still, in such cases, we show that a weak Erlang memory leads to convergence to an endemic equilibrium, while a strong Erlang memory can determine the appearance of oscillations.

Our work has, of course, a number of limitations. The qualitative analysis is limited to local stability, while global stability results are still an open problem. Another limitation is that we have considered only Erlang memory kernels, while non-Erlang distributions could be required to represent the lower time-scale of the process of formation and acquisition of information, as compared to the memory fading of acquired information \cite{DONOFRIO2022}. We also remark that in this paper, we extend minimal behavioural epidemic models that, while interesting as prototypes of phenomena that are not captured by the corresponding non-behavioural analogues, are  unable to describe complex dynamics such as those produced, for example, by the spread of SARS-CoV-2. Future work will be devoted to the extension of more complex models: we plan to consider also the case of models, including non-mandatory pharmaceutical interventions (e.g. vaccination and treatments) whose adoption by individuals is usually strongly affected by opinion-driven change of human behaviour.
Finally, the impact of host heterogeneity on epidemic dynamics has long been acknowledged \cite{brauer2008epidemic,diekmann1990definition}. Recent literature introduces general Kermack and McKendrick's epidemic models that consider heterogeneity, for a population characterized by different traits, with different activity levels or disease parameters \cite{bootsma2023effect,bootsma2023separable,brauer2009age}. The provided formula expresses the FoI, delineating contributions from individuals infected at a specific time and with a particular trait. This suggests, as a natural extension, a generalization of the constitutive law for the FoI with an information index aimed at exploring how heterogeneity influences epidemic aspects in a model that incorporates human behaviour and broadens the analysis presented in this paper.

\bigskip

\paragraph*{Acknowledgements}
The authors sincerely thanks Alberto d'Onofrio (University of Trieste) for his valuable suggestions and discussions.
This work has been performed under the auspices of the Italian National Group for Scientific Computing (GNCS) and  the
Italian National Group for Mathematical Physics (GNFM) of the
National Institute for Advanced Mathematics (INdAM).\\ 
B.B. acknowledges EU funding within the NextGenerationEU---MUR PNRR Extended Partnership initiative on Emerging Infectious Diseases (Project no. PE00000007, INF-ACT).\\
B.B. also acknowledges PRIN 2020 project (No. 2020JLWP23) ``Integrated Mathematical Approaches to Socio--Epidemiological Dynamics''. \\

\paragraph*{Data availability statement} Data sharing not applicable to this article as no datasets were generated or
analysed during the current study.

\appendix

\section{Discretization scheme}
\label{AppendixA}

Consider an uniform mesh $t_n=nh$, where $n=0,1,\dots$, and $h>0$ is the stepsize. We define the following discretization scheme for the continuous model \eqref{mod continuo22}:

\begin{align}\label{mod discreto}
	S_{n+1} &= S_n +h \big( \lambda - \mu S_{n+1}-\beta(M_n) S_{n+1}F_n \big), \nonumber \\
	F_{n+1} &= F_0(t_{n+1})+h \sum_{j=0}^n  A_{\mu}(t_{n+1-j})\beta(M_j)F_jS_{j+1}, \\
	M_{n+1} &= M_0(t_{n+1})+h \sum_{j=0}^n  K(t_{n+1-j})g(F_j), \nonumber
\end{align}
for $ n=0,1,\dots$, where the forcing functions $F_0(t)$ and $M_0(t)$ are given by \eqref{ForcingF00} and $\eqref{forcingFM}_2$, respectively.  Here, $S_n \approx S(t_n)$, $F_n \approx F(t_n)$, $M_n \approx M(t_n)$ for $ n=0,1,\dots$ and   $S_0=S(0)$, $F_0=F_0(0)$ and $M_0=M_0(0)$. This numerical method is called non-standard because it is based on a modified rectangular rule, which is implicit in $S$ and explicit in $F$ and $M$.
A detailed analysis of the qualitative and asymptotic behaviour of the numerical solution to \eqref{mod continuo22}, obtained by \eqref{mod discreto}, is presented in \cite{DISCRETONOSTRO}. Here we enunciate the main results.

Denote by 
$$\textbf{\textit{e}}(h;t_n)= \begin{bmatrix}
	S(t_n) \\
	F(t_n) \\
	M(t_n)
\end{bmatrix}- \begin{bmatrix}
	S_n \\
	F_n \\
	M_n
\end{bmatrix}, $$
the global error of discretization \eqref{mod discreto}. 
%Consider the interval $[0,T]$ and choose $h=T/M$, where $M$ is a positive integer. 

\begin{theorem}
	Assume that the infectivity function $A(\tau)$  is continuously differentiable on an interval $[0,T]$ and that $\{S_n\}_{n\in\mathbb{N}_0}$, $\{F_n\}_{n\in\mathbb{N}_0}$ and $\{M_n\}_{n\in\mathbb{N}_0}$ are the approximations to the solution of \eqref{mod continuo22}, defined by \eqref{mod discreto}. Let $h=T/\Bar{n}$, where $\Bar{n}$ is a positive integer, then
 
	$$
	\lim_{h \rightarrow 0} \max_{0 \leq n \leq \Bar{n}} ||\textbf{e}(h; t_n)|| =0.
	$$
 
\noindent Furthermore, the convergence rate is linear.
\end{theorem}

\begin{theorem}
	Consider equation \eqref{mod discreto}, the following inequalities hold:  $0<S_n\leq S_{\max}(h),$ $0\leq F_n \leq F_{\max}(h)$ and $0 \leq M_n \leq M_{\max}(h)$, for $n =0,1,\dots$ and $h>0$. Here $S_{\max}(h),F_{\max}(h)$ and $M_{\max}(h)$ are positive constants.
 %Let $(S_n, F_n, M_n)$ be the solution to the discrete equation \eqref{mod discreto}, with $h > 0$ and with the initial values $S_0, F_0, M_0>0$, then $S_n, F_n, M_n$  are non-negative for all $ n \geq 0$. 
\end{theorem}

\noindent The definition of the numerical basic reproduction number

$$R_0(h)= \frac{\lambda}{\mu} \Bar{A}_{\mu}(1,h),$$ 

\noindent with $\Bar{A}_{\mu}(1,h)=h \sum_{j=1}^{\infty} A_{\mu}(t_{j})$, as the discrete equivalent to \eqref{R0}, allows an investigation about the equilibria of the discrete system \eqref{mod discreto}, which parallels the one reported in Sections \ref{sec4}-\ref{sec5}. In particular, if $R_0(h)>1$, then  system \eqref{mod discreto} admits a unique endemic equilibrium that we denote by $EE(h)=(S_{e}(h),F_{e}(h), M_{e}(h))$.  The endemic equilibria are obtained from the following conditions:

\begin{equation*}
     M_{e}(h)=\Bar{K}(1,h)g(F_{e}(h)), \quad \quad S_{e}(h)= \frac{\lambda}{\mu R_0(h) \beta(\Bar{K}(1,h)g(F_{e}(h))},
\end{equation*}

\noindent and $F_{e}(h)$ is the solution of the nonlinear equation

\begin{equation*}
  x-\mu R_0(h)+\frac{\mu}{\beta(\Bar{K}(1,h)g(x))}=0,
\end{equation*}
\noindent  in the interval $[0,F_{max}(h)]$, where $\Bar{K}(1,h)=h \sum_{j=1}^{\infty} K(t_{j})$. Set $\Bar{A}_{\mu}(z,h)=h \sum_{j=1}^{\infty} A_{\mu}(t_{j}) z^j$ and $\Bar{K}(z,h)=h \sum_{j=1}^{\infty} K(t_{j}) z^j$, the following theorem gives a sufficient condition for the local asymptotic stability of the equilibria.

\begin{theorem}\label{th Stab discr}
	Let be $h>0$. Any equilibrium $EE(h)$ of model \eqref{mod discreto} is locally asymptotically stable if  the characteristic equation 

 {\small \begin{equation*}\label{EEequa1app}
		\beta(M_{e}(h)) S_{e}(h) \Bar{A}_{\mu}(z,h) \bigg( 1+ \frac{g'(F_{e}(h)) \beta'(M_{e}(h))F_{e}(h)\Bar{K}(z,h)}{\beta(M_{e}(h))} \bigg) = 1+ h \frac{\beta(M_{e}(h))F_{e}(h)} {1-z + \mu h},
	\end{equation*}}
 
\noindent has no solution for $|z| \leq 1$. 
%For $|z|=1$, \eqref{EEequa1app} has be interpreted as the limit for $\Tilde{z}\rightarrow z$, $|\Tilde{z}|<1$.
\end{theorem}

\noindent From Theorem \ref{th Stab discr}  sufficient conditions, written in terms of the information parameters and functions, for the numerical endemic equilibrium to be locally stable, are derived. The reader is referred to \cite{Nostro}. A comparison with the results of Section \ref{sec5} confirms that the numerical solution replicates the stability properties of the continuous one.
 Furthermore, as $h \rightarrow 0$ the accuracy with which the numerical method approximates the entire dynamics of the continuous system \eqref{mod continuo22} increases.

%\begin{theorem}
%Assume that:
    %\begin{itemize}
       % \item  $K(\tau)=\delta(\tau)$ and $A_{\mu}(\tau)$ is positive definite, then  $EE(h)$ is LAS for any choice of the information functions $\beta$ and $g$;

      %  \item  $K(\tau)=\delta(\tau)$ or $K(\tau)=a e^{-a \tau}$ and that functions $\beta$ and $g$ satisfy the inequality

       % \begin{equation*}  g'(F_{e}(h))\beta'(M_{e}(h))F_{e}(h) \geq -2 \beta(M_{e}(h)) \phi(h),
    % \end{equation*}

       % where $\phi(h)=\big( \frac{e^{ah}-1}{ah} \big) \frac{1}{e^{ah}}$. Then, for $h>0$,  $EE(h)$ is LAS for any choice of $A_{\mu}(\tau)$.
        
   % \end{itemize}
%\end{theorem}

\section{Paley-Wiener Theorem}
\label{AppendixB}
The theory on the asymptotic stability of Volterra integral  equations of convolution type is essentially based on the Paley-Wiener theorem \cite[p. 59]{PW1934}. An extension of this result to the case where the matrix kernel does not belong to $L^1(0,\infty)$ is provided by Ch. Lubich \cite[Th. 9.1]{LUBICH1983}. Here, we report this result for the sake of completeness.

\begin{theorem}\label{pali1}  Consider the system of Volterra integral equations 
\begin{equation}\label{lineare generale}
     \textbf{y}(t)=\textbf{f}(t)+ \int_0^t \mathcal{\textbf{Q}}(t-\tau) \textbf{y}(\tau) ds,
 \end{equation}
where the matrix $\mathcal{\textbf{Q}}(\cdot)$ and the function $\textbf{f}(t)$ satisfy

\begin{equation}\label{hp lineare gen}
    \mathcal{\textbf{Q}}(\cdot)= \textbf{B}(\cdot) + \textbf{A}_{\infty}, \quad \textbf{A}_{\infty}=\underbrace{\begin{bmatrix}
\textbf{A}^1_{\infty} \\
0
\end{bmatrix}}_{d}\hspace{-.3cm}\begin{array}{l}
\} l \\
\} k
\end{array}, \quad \textbf{f}(t)= \begin{bmatrix}
\textbf{f}_0(t) \\
\textbf{f}_1(t)
\end{bmatrix}\hspace{-.3cm}\begin{array}{l}
\} l \\
\} k
\end{array},
\end{equation}
with $\textbf{B}(\cdot) \in L^1(0, \infty)$ and $k=d-l=$dimension of the null space of $\textbf{A}_{\infty}$. Then $\textbf{y}(t) \rightarrow 0$, whenever $\textbf{f}^{\prime}_0(t),$ $\textbf{f}_1(t) \longrightarrow 0$ for $t\rightarrow \infty$, if and only if the equation
\begin{equation}\label{det}
    w^l \text{det} \Big( I - \int_0^\infty \mathcal{\textbf{Q}}(t) e^{-w t} dt \Big) = 0, 
\end{equation}
has no solution for $ Re(w) \geq 0$. 
\end{theorem}
\begin{remark}
    For $Re(w)=0$, equation \eqref{det} must be interpreted as limit for $\Tilde{w} \rightarrow w$, $Re(\Tilde{w}) >0$.
\end{remark}

\bibliographystyle{abbrv}
\bibliography{bibliog}

%\end{linenumbers}

\end{document}